\documentclass[12pt]{amsart}

\usepackage[USenglish]{babel}
\usepackage{amsmath,amsthm,amssymb,amsfonts}
\usepackage{mathtools}
\usepackage{mathrsfs}
\usepackage{bbm}

\usepackage{indentfirst}
\usepackage{microtype}
\usepackage{enumerate}
\usepackage{accents}
\usepackage{cancel} 

\usepackage{multirow}
\usepackage{longtable}

\usepackage{float}
\usepackage{graphicx}
\usepackage{placeins}
\usepackage{caption}

\usepackage{tikz}
\usetikzlibrary{patterns}
\usepackage{tikz-cd}

\usepackage{hyperref}

\allowdisplaybreaks

\newtheorem{thm}{Theorem}[section]

\newtheorem{lem}[thm]{Lemma}
\newtheorem{prop}[thm]{Proposition}
\newtheorem{conj}[thm]{Conjecture}

\theoremstyle{definition}

\theoremstyle{remark}

\newtheorem*{ntt}{Notation}

\numberwithin{equation}{section}

\renewcommand{\P}{\mathbb{P}}      

\newcommand{\eps}{\varepsilon}   

\newcommand\restr[2]{{           
  \left.\kern-\nulldelimiterspace #1%
  \right|_{#2}%
 }}

\newcommand\minus{               
  \setbox0=\hbox{-}%
  \vcenter{%
    \hrule width\wd0 height \the\fontdimen8\textfont3%
  }%
}

\newcommand{\Var}{\mathrm{Var}}  
\newcommand{\h}{\mathrm{h}}      

\newcommand{\gtil}{\widetilde{g}}

\tikzset{                        
    symbol/.style={%
        ,draw=none
        ,every to/.append style={%
            edge node={node [sloped, allow upside down, auto=false]{$#1$}}}
    }
}

\newcommand{\uptxt}[2]{          
  \mathrel{\overset{\makebox[0pt]{\mbox{\normalfont\tiny #1}}}{#2}}}

\newcommand{\negphantom}[1]{\settowidth{\dimen0}{#1}\hspace*{-\dimen0}}

\restylefloat{table}

\begin{document}

\title[Elementary heuristic for extended Goldbach]{An elementary heuristic for Hardy--Littlewood extended Goldbach's conjecture}%
\author{Christian T\'{a}fula}%
\address{RIMS, Kyoto University, 606-8502 Kyoto, Japan}%
\email{tafula@kurims.kyoto-u.ac.jp}%



\begin{abstract}
 The goal of this paper is to describe an elementary combinatorial heuristic that predicts Hardy \& Littlewood's extended Goldbach's conjecture. We examine common features of other heuristics in additive prime number theory, such as Cram\'{e}r's model and density-type arguments, both of which our heuristic draws from. Apart from the prime number theorem, our argument is entirely elementary, in the sense of not involving complex analysis. The idea is to model sums of two primes by a hypergeometric probability distribution, and then draw heuristic conclusions from its concentration behavior, which follows from Hoeffding-type bounds.
\end{abstract}
\maketitle

\section{Introduction}
 Denote by $\mathbb{N}$ the set of natural numbers with $0$. Writing $\P\subseteq\mathbb{N}$ for the set of prime numbers, given $h\geq 1$, we consider its \emph{representation functions}:
 \begin{align*}
  r_{\P,h}(n) &:= \{(k_1,\ldots,k_h) \in \P^h : \textstyle{\sum}_{i=1}^{h} k_i = n \} &&\text{ for } n\in\mathbb{N}, \\
  s_{\P,h}(x) &:= \{(k_1,\ldots,k_h) \in \P^h : \textstyle{\sum}_{i=1}^{h} k_i \leq x \} = \sum_{n\leq x} r_{\P,h}(n) &&\text{ for }x \in\mathbb{R}_{\geq 0}.
 \end{align*}
 These quantities are counting the number of solutions to $k_1+\ldots+k_h=n$ and $k_1+\ldots+k_h \leq x$, respectively, with $k_i \in \P$ for $i=1,\ldots, h$, considering repetitions. As usual, we denote the \emph{prime counting function} by $\pi(x):= s_{\P,1}(x)$, and also the \emph{characteristic function} of $\P$ by $\mathbbm{1}_{\P}(n) := r_{\P,1}(n)$. These functions are generated by the following power series:
 \begin{equation*}
  \left(\sum_{p\in\P} z^{p} \right)^h = \sum_{n\geq 0} r_{\P,h}(n)z^n, \qquad \frac{\left(\sum_{p\in\P}z^{p} \right)^h}{1-z} = \sum_{n\geq 0} s_{\P,h}(n)z^n.
 \end{equation*}
 Notice that this same setup can be replicated for any set $\mathscr{A}\subseteq\mathbb{N}$ in place of $\P$; in this paper, however, we focus on the set of primes.
 
 From Chebyshev's elementary estimate $\pi(x) = \Theta(x\log(x)^{-1})$, one can, with relative ease, show that $s_{\P,h}(n) = \Theta(n^h\log(n)^{-h})$ for all $h\geq 1$. From this, it seems reasonable, therefore, to expect that the growth rate of $r_{\P,h}(n)$ should approximate $n^{h-1}\log(n)^{-h}$ as $n$ gets large, i.e. similar to $s_{\P,h}(n)/n$. Indeed, in Hardy \& Littlewood's 1923 seminal paper \emph{``Some problems on 'Partitio Numerorum' III''} \cite{harlit23} it is shown, assuming a weak version of the Generalized Riemann Hypothesis (GRH), that for all $h\geq 3$, when $n\to +\infty$ through integers with same parity than $h$, the following holds:
 \begin{equation}
  r_{\P,h}(n) \sim \frac{2C_h}{(h-1)!}\frac{n^{h-1}}{\log(n)^h}\prod_{\substack{p\mid n \\ p\geq 3}} \left(\frac{(p-1)^h+(-1)^h(p-1)}{(p-1)^h -(-1)^h} \right), \label{hdlw}
 \end{equation}
 where $C_h$ is the constant
 \begin{equation*}
  C_h := \prod_{p\geq 3} \left(1- \frac{(-1)^h}{(p-1)^h}\right).
 \end{equation*}
 To be precise, their assumption was that there is some $\eps>0$ such that, for every Dirichlet character $\chi$, it holds:
 \begin{equation}
  \text{``If } L(s,\chi)= 0, \text{ then } \mathrm{Re}(s)\leq 3/4-\eps.\text{''} \label{wgrh}
 \end{equation}
 We are not going to describe this in detail, for $L$-functions will not be essential to our discussion. Nevertheless, roughly 15 years later, I. M. Vinogradov \cite{vinogradov54} introduced a powerful new technique to prove \eqref{hdlw} unconditionally. This was a remarkable achievement in additive prime number theory, implying, in particular, that every large even integer is the sum of at most four primes. Details on Vinogradov's method applied to the even more general \emph{Waring-Goldbach's problem} may be found in Hua \cite{hua65}. It is interesting to remark that Vinogradov's method for $h=3$ can also be used to show that ``almost all'' even numbers may be written as a sum of two primes. More precisely, in 1975, Montgomery \& Vaughan \cite{monvau75} were able to ensure the existence of an effectively computable small $\delta>0$ for which $|\mathfrak{E}\cap[0,x]| \ll x^{1-\delta}$, where $\mathfrak{E}:= \{2n:n\in\mathbb{N}\}\setminus (\P+\P)$\footnote{In general, given $\mathscr{A},\mathscr{B}\subseteq\mathbb{N}$ we have the \emph{sumset} $\mathscr{A}+\mathscr{B} = \{a+b:a\in\mathscr{A},b\in\mathscr{B}\}$, and for $h\geq 2$ the \emph{$h$-fold sumset} of $\mathscr{A}$ is written as $h\mathscr{A} := \underbrace{\mathscr{A} + \ldots + \mathscr{A}}_{h \text{ times}}$.} is the \emph{exceptional Goldbach's set}. The full Goldbach's conjecture consists of the far stronger statement:
 
 \theoremstyle{plain}
 \newtheorem*{gold}{Goldbach's Conjecture}
 \begin{gold}\label{goldbach}
  $\mathfrak{E} = \{0,2\}$.
 \end{gold}

 In spite of that, notice that the estimate in \eqref{hdlw} is only proven for $h > 2$. When it comes to $h=2$, even assuming the full GRH (with ``$1/2$'' instead of ``$3/4-\eps$'' at \eqref{wgrh}), Hardy \& Littlewood \cite{harlit24} were only able to show that $|\mathfrak{E}\cap[1,x]|\ll_\eps x^{1/2 +\eps}$. Based on that, the following is then stated as a conjecture:

 \begin{conj}[Hardy \& Littlewood]\label{aconj}
  When $n\to +\infty$ through the even numbers,
  \begin{equation}
   r_{\P,2}(n) \sim 2C_2\frac{n}{\log(n)^2}\prod_{\substack{p\mid n \\ p\geq 3}} \left(\frac{p-1}{p-2}\right), \label{hlgbcj}
  \end{equation}
  where $C_2$ is the constant
  \begin{equation}
   C_2 := \prod_{p\geq 3} \left(1- \frac{1}{(p-1)^2}\right) \approx 0.6601618\ldots \label{c2}
  \end{equation}
 \end{conj}

 This would imply, in particular, that $\mathfrak{E}$ is finite. Our goal in this paper is to motivate this conjecture with far less machinery; in particular, avoiding any non-elementary sieve-theoretic estimates or Fourier analysis. Apart from the Prime Number Theorem (PNT), our argument is completely elementary, in the sense that it avoids complex analysis. On the probabilistic side, we only use basic facts about the hypergeometric distribution, together with a Hoeffding-type inequality due to V. Chv\'{a}tal \cite{chv79}. For the number-theoretic estimates used, refer to Hardy \& Wright \cite{hardy08}.
 
 \begin{ntt}
  We use the standard asymptotic symbols $\Theta,\asymp,O,\ll,o,\sim$, as well as the ``floor'' and ``ceiling'' functions $\lfloor\cdot\rfloor$, $\lceil \cdot\rceil$. The symbol ``$\approx$'' is used to denote rough, conjectural or heuristic approximations, the latter being indicated by a question mark ``$?$'' above it.
  
  We use $\Pr$ for the probability measure, $\mathbb{E}$ for expectation and $\Var$ for variance. Given a probability space $(\Omega, \mathcal{F}, \Pr)$, a random variable (abbreviated r.v.) $X:\Omega\to \mathbb{R}$ is said to be a \emph{Bernoulli trial} when $X(\Omega) \subseteq \{0,1\}$. The conditional probability of an event $E$ given another event $F\in\mathcal{F}$ is denoted by $\Pr(E|F) := \Pr(F)^{-1}\cdot \Pr(E\wedge F)$. Whenever we write an asymptotic sign with a superscripted ``a.s.'', we mean that the corresponding limit holds \emph{almost surely}, i.e. for a subset of $\Omega$ with complement having measure $0$.
  
  In general, letters $p$, $q$ denote primes, and $n\mid m$ means ``$n$ divides $m$''. We denote by $(a,b)$ the $\gcd$ of $a$ and $b$, $\pi(x) := |\P\cap[0,x]|$ the prime counting function, $\varphi(n):= |\{1\leq k\leq n : (n,k)=1\}|$ is Euler's totient function, and $\omega(n) := |\{p\in\P: p\mid n\}|$ the distinct prime factors counting function.
 \end{ntt}
 
\section{Flaws in simpler heuristics}
 Before presenting our argument, it will be instructive to look at two related heuristic arguments that happen to fail in predicting Hardy \& Littlewood's expected estimate in Conjecture \ref{aconj}. The purpose is to give an overview of the main ingredients involved in heuristics concerning prime numbers, and also to contextualize our model with related concepts in the literature.
 
 Every heuristic argument can be said to have some sort of \emph{leap of faith}, with some more sound than others. We draw attention to this keyword, for we shall use it to describe the non-rigorous step involved in heuristic arguments in general. Such step, usually, is backed by some probabilistic reasoning or density-type argument which leads to an apparently sensible conclusion. Perhaps the most influential heuristic argument in analytic number theory is the one formulated by H. Cram\'{e}r \cite{cra36}, which we briefly outline.

\subsection{Na\"{i}ve Cram\'{e}r's model}
 Commonly referred in the literature as \emph{Cram\'{e}r's model}, this heuristic argument consists of considering a random subset $\mathcal{R}\subseteq\mathbb{N}$ with $\Pr(n\in\mathcal{R}) = \log(n)^{-1}$, the events ``$n\in\mathcal{R}$'' being mutually independent. The idea is, then, to infer non-multiplicative properties of the primes based on probabilistic aspects of $\mathcal{R}$. The $\log(n)^{-1}$ term comes from the PNT, which states that $\pi(x)/x \sim \log(x)^{-1}$.\footnote{cf. Theorem 6, p. 10 of Hardy \& Wright \cite{hardy08}.} Among its many variations, the simplest one, which was first described by Cram\'{e}r, is sometimes referred to as the \emph{na\"{i}ve Cram\'{e}r's model}, and it is the one we are going to be focusing on.
 
 This concept is not as unsound as it may seem at first glance. Similar constructions are commonplace in probabilistic combinatorics (cf. \cite{deshenlan98, erdo56, erdren60, erdtet90, lan95, vvu00wp}), and a classical exposition of related concepts may be found in Chapter III of Halberstam \& Roth \cite{halberstam83}. To be a bit more rigorous, the general recipe of such constructions can be described as follows. For any sequence of real numbers $(\alpha_n)_{n\geq 0}$ with $\alpha_n \in [0,1]$, it is possible to construct a probability space $(2^{\mathbb{N}},\mathcal{F},\Pr)$ satisfying:
 \begin{enumerate}[(i)]
  \item The events $E_n:=\{\mathcal{R}\subseteq\mathbb{N}: n\in\mathcal{R}\}$ are measurable and $\Pr(E_n) = \alpha_n$;

  \item $\{E_0, E_1, E_2,\ldots\}$ is a collection of mutually independent events;

  \item $\mathcal{F}$ is the $\sigma$-algebra induced by the collection $\{E_0, E_1, E_2, \ldots\}$.
 \end{enumerate} 
 The well-definedness of this space follows once infinite Cartesian products of probability spaces are established,\footnote{This may be found in several sources, most notably in Halmos \cite{halmos74} (cf. Section 38). In Section III of Halberstam \& Roth \cite{halberstam83} one finds this construction for the specific case we are working with.} for we may also interpret it as the product space of the Bernoulli trials $\mathbbm{1}_{\mathcal{R}}(n)$ with $\Pr(\mathbbm{1}_{\mathcal{R}}(n)=1) = \alpha_n$. It is also common to require the series $\sum_{n\geq 0} \alpha_n$ to diverge so as to work solely with infinite subsets of $\mathbb{N}$, for, in this case, $\mathcal{R}$ will almost surely be infinite. This follows from the Borel--Cantelli lemma, which is a fundamental tool in the subject.
 
 \theoremstyle{plain}
 \newtheorem*{bc}{Borel--Cantelli lemma}
 \begin{bc}[Section VIII.3, p. 200 of Feller \cite{feller1}]
  Let $(E_n)_{n\geq 0}$ be a sequence of events in a probability space. The following holds:
  \begin{enumerate}[(i)]
   \item If $\sum_{n=0}^{\infty} \Pr(E_n) < +\infty$, then, with probability $1$, only a finite number of these events take place;
    
   \item If $\sum_{n=0}^{\infty} \Pr(E_n)$ diverges and the events are independent, then, with probability $1$, an infinite number of these events occur.
  \end{enumerate}  
 \end{bc}
 
 The idea behind these probability measures is to study sets with certain prescribed rates of growth. Indeed, with a specific version of the strong law of large numbers\footnote{See Theorem 3.1 in \cite{taf19} for a short proof of this version of the strong law.} one deduces that, when $\sum_{n\geq 0} \alpha_n$ diverges,
 \[ |\mathcal{R}\cap [0,x]| \uptxt{a.s.}{\sim} \sum_{n\leq x} \alpha_n; \]
 that is, this asymptotic relation will hold for \emph{almost all} subsets $\mathcal{R}\subseteq \mathbb{N}$, i.e. all but a set with measure $0$. Under this framework, the so-called \emph{na\"{i}ve} Cram\'{e}r's model is just a version of this construction done with the intent of having $|\mathcal{R}\cap[0,x]| \uptxt{a.s.}{\sim} x\log(x)^{-1}$. In 1936, Cram\'{e}r \cite{cra36} used his model to conjecture that gaps between primes must be mostly small; more precisely, that $p_{n+1}-p_n = O(\log(p_n)^2)$, where $p_n$ is the $n$-th prime. He also showed that, assuming the Riemann Hypothesis, it is possible to derive $p_{n+1}-p_n = O(\sqrt{p_n}\log(p_n))$, which is still much weaker than the former, for which all current numerical data suggests to be really the case (cf. Section A.8 of Guy \cite{guy94}).
 
 One way to explore this model in Goldbach's problem is by considering the r.v. $r_{\mathcal{R},2}(n) = \sum_{k\leq n} \mathbbm{1}_{\mathcal{R}}(k)\mathbbm{1}_{\mathcal{R}}(n-k)$ and studying its distribution. A first step is calculating the expected value $\mathbb{E}(r_{\mathcal{R},2}(n))$. Considering only integers greater than $3$ in $\mathcal{R}$, we have, by independence:
 \begin{equation*}
  \mathbb{E}(r_{\mathcal{R},2}(n)) = \sum_{3\leq k\leq n-3} \frac{1}{\log(k)}\frac{1}{\log(n-k)} \sim 2\underbrace{\int_{3}^{n/2} \frac{\mathrm{d}t}{\log(t)\log(n-t)}}_{=:\,I_n}.
 \end{equation*}
 Moreover, since
 \[ \frac{2}{\log(n-3)}\int_{3}^{n/2} \frac{\mathrm{d}t}{\log(t)} \leq 2I_n \leq \frac{2}{\log(n)-\log(2)}\int_{3}^{n/2} \frac{\mathrm{d}t}{\log(n-t)}, \]
 we may deduce, in view of $\int_{3}^{x} \log(t)^{-1}\mathrm{d}t \sim x\log(x)^{-1}$, that
 \begin{equation*}
  \mathbb{E}(r_{\mathcal{R},2}(n)) \sim \frac{n}{\log(n)^2}.
 \end{equation*}
 With the aid of some concentration inequalities such as Chernoff bounds (cf. Theorem 1.8, p. 11 of Tao \& Vu \cite{tao06}) it is even possible to show that $r_{\mathcal{R},2}(n)\uptxt{a.s.}{\sim} n\log(n)^{-2}$. The \emph{leap of faith} here would then be to assume that the sequence of primes is an ``average sequence'' in this space, and hence 
 \begin{equation}
  r_{\P,2}(n) \uptxt{?}{\approx} \frac{n}{\log(n)^{2}}. \label{wrng1}
 \end{equation}
 
 Substituting ``$\approx$'' by ``$\sim$'', however, not only produces the wrong\footnote{That is, assuming Hardy-Littlewood's conjecture do hold.} constant, but the wrong growth order! One of the reasons for this is that this model does not take into account the distribution of primes among residue classes. It leads us to expect, for example, infinitely many ``conjoined twin primes'' $(p, p+1)$, which is far from actuality. A more in-depth analysis on further discrepancies produced by this model may be found in Pintz \cite{pin07}. A modified version of Cram\'{e}r's model which does take divisibility by small primes into account can be seen in Granville \cite{gra95}. Further refinements to this probabilistic reasoning have been extensively explored by other authors, the most prominent example being the far-reaching conjecture concerning the distribution of prime numbers over systems of integer-valued polynomials by P. T. Bateman \& R. A. Horn \cite{bathor62}.
 
\subsection{A sieve theory inspired attempt}
 In Section 4 of \cite{harlit23}, Hardy \& Littlewood compare their conjectured estimate with other heuristic arguments available at the time that arrive at a different conclusion, in particular, one derived from sieve theory which is generally attributed to V. Brun.\footnote{Although following from Brun's pure sieve, the heuristic concerning Goldbach was not formulated by Brun, as remarked in footnote 5 at p. 33 of Hardy \& Littlewood \cite{harlit23}.} We will not describe his method, but we will present a rough sketch of the idea behind it. First, by using only Chebyshev's elementary PNT\footnote{cf. Theorem 4.6, p. 82 of Apostol \cite{apostol76}.} $\pi(x) =\Theta(x \log(x)^{-1})$, we derive
 \begin{align}
  r_{\P,2}(n) &= \sum_{k\leq n} \mathbbm{1}_{\P}(k)\mathbbm{1}_{\P}(n-k) \nonumber \\
  &= \sum_{\sqrt{n}<k\leq n-\sqrt{n}} \mathbbm{1}_{\P}(k)\mathbbm{1}_{\P}(n-k) + O\left( \frac{\sqrt{n}}{\log(n)}\right). \label{hgold1}
 \end{align}
 Write ``$x \stackrel{m}{\equiv} a$'' as short for ``$x\equiv a~ \mathrm{(mod}~m\mathrm{)}$'', and fix $n\geq 8$ an even integer. When $\sqrt{n}<k\leq n-\sqrt{n}$, it holds:
 \begin{align*}
  \mathbbm{1}_{\P}(k)\mathbbm{1}_{\P}(n-k) = 1 \iff \forall p\leq \sqrt{n}, ~k(n-k) \stackrel{p}{\not\equiv} 0 \iff \forall p\leq \sqrt{n}, ~\begin{cases} k \stackrel{p}{\not\equiv} 0 \text{ and} \\ k \stackrel{p}{\not\equiv} n.\end{cases}
 \end{align*}
 Hence, in order to solve Goldbach's conjecture, all we would need to do would be to \emph{sift out} the undesirable $k$s from $(\sqrt{n}, n-\sqrt{n}]$. More precisely, we need to estimate
 \begin{equation}
  \mathfrak{G}(n) := \left|\left\{k\leq n: k \stackrel{p}{\not\equiv} 0 \text{ and } k \stackrel{p}{\not\equiv} n,~\forall p\leq \sqrt{n} \right\}\right|, \label{hgold3}
 \end{equation}
 which is equal to the main term in \eqref{hgold1}.
 
 Writing $\mu$ for the \emph{M\"{o}bius function}:
 \begin{equation*}
  \mu(n) := \begin{cases}
              (-1)^{\omega(n)} \text{ if } n \text{ squarefree,} \\
              0 \phantom{(-1)^{\omega(n)}}\negphantom{$0$}\text{ otherwise;}
            \end{cases}
 \end{equation*}
 one could effectively estimate $\mathfrak{G}(n)$ in terms of $\mu$ using the identity
 \[ |\{k\leq x : (k,n) = 1\}| = \sum_{d\mid n} \mu(d)\left\lfloor \frac{x}{d}\right\rfloor, \]
 which is basically a fancy version of Eratosthenes' sieve. Both Brun's and Selberg's sieve stem from a careful analysis of the error term produced by substituting ``$\lfloor x/d\rfloor$'' by ``$x/d$''. A detailed overview of both methods may be found in Chapter IV of Halberstam \& Roth \cite{halberstam83}. As sieve theory per se is not the main goal of this exposition, we shall consider a more picturesque approach.
 
 Writing $x\# := \prod_{p\leq x} p$ for the \emph{primorial} function, we have, by the Chinese remainder theorem,
 \begin{equation*}
  \left|\left\{k\leq \sqrt{n}\#: k \stackrel{p}{\not\equiv} 0 \text{ and } k \stackrel{p}{\not\equiv} n,~\forall p\leq \sqrt{n} \right\}\right| = \prod_{p\mid n} (p-1) \prod_{\substack{p\,\nmid\, n \\ p\leq \sqrt{n}}}(p-2);
 \end{equation*}
 thus, a \emph{leap of faith} in this direction should look something like:
 \begin{align}
  \mathfrak{G}(n) &\stackrel{?}{\approx} \frac{n}{\sqrt{n}\#} \prod_{p\mid n} (p-1) \prod_{\substack{p\,\nmid\, n \\ p\leq \sqrt{n}}}(p-2) \label{hgold5} \\
  &= n \prod_{p\mid n} \left(1-\frac{1}{p}\right) \prod_{\substack{p\,\nmid\, n \\ p\leq \sqrt{n}}}\left(1-\frac{2}{p} \right). \nonumber
 \end{align}
 Put into words, we are assuming that the distribution among residues provided by the Chinese remainder theorem is equidistributed on very short intervals. Moreover, since it is possible to deduce from Mertens' 3rd theorem\footnote{cf. Theorem 429, p. 466 of Hardy \& Wright \cite{hardy08}.} that
 \begin{equation*}
  \prod_{3\leq p\leq \sqrt{n}}\left(1-\frac{2}{p} \right) = \prod_{3\leq p\leq \sqrt{n}}\left(1-\frac{1}{p} \right)^2 \left(1-\frac{1}{(p-1)^2} \right) \sim \frac{16e^{-2\gamma}}{\log(n)^2}C_2, \label{hgold52}
 \end{equation*}
 where $\gamma:= \lim_n (\sum_{k=1}^n 1/k - \log(n)) \approx 0.57721\ldots$ is the \emph{Euler-Mascheroni constant} and $C_2$ is the constant from \eqref{c2}, we may rewrite \eqref{hgold5} in terms of $r_{\P,2}(n)$ as follows:
 \begin{equation}
  r_{\P,2}(n) \stackrel{?}{\approx} 8e^{-2\gamma}C_2\frac{n}{\log(n)^2} \prod_{\substack{p\mid n \\ p\geq 3}} \left(\frac{p-1}{p-2}\right). \label{hgold6}
 \end{equation}

 In contrast to the na\"{i}ve Cramer's model, this heuristic at least agrees with Hardy \& Littlewood's conjectured growth order, missing its constant only by a factor of $4e^{-2\gamma}$. That is exactly why substituting ``$\approx$'' by ``$\sim$'' in \eqref{hgold6} \emph{must} be unsound! In Section 4 of Hardy \& Littlewood \cite{harlit23} it is shown that if
 \[ r_{\P,2}(n) \sim K\frac{n}{\log(n)^2} \prod_{\substack{p\mid n \\ p\geq 3}} \left(\frac{p-1}{p-2}\right) \]
 for some real constant $K>0$, then $K$ \emph{must} be equal to $2C_2$. A possible explanation for the unsoundness of this substitution could be that the regularity we are requiring in \eqref{hgold5} is for an, indeed, \emph{very} short interval, as in logarithmically short, for it follows from the PNT that $\sqrt{n}\# \sim e^{(1+o(1))\sqrt{n}}$. Notwithstanding, this is a more prolific approach than the previous one. Roughly a decade before Vinogradov's theorem on the sum of three primes, L. G. Schnirelmann proved using purely combinatorial methods that every integer $n\geq 2$ can be written as sum of at most $K$ primes, where $K>0$ is some large, effectively computable constant not depending on $n$. This was the first significant result on the direction of Goldbach's conjecture, and at the heart of his proof is the fact that substituting ``$\approx$'' by ``$\ll$'' in \eqref{hgold5} is a sound substitution (cf. Theorem 7.2, p. 186 of Nathanson \cite{nathanson96}).

\section{The Hypergeometric Model}
 We now describe the framework of our heuristic. Our argument is a sort of hybrid between the two we have just presented, dealing probabilistically with $p\mid 2n$ and making a density-type argument for $p\nmid 2n$, culminating, first, in \eqref{mainleap}. The \emph{leap of faith}, then, will come from the concentration behavior of the r.v.s used to model $r_{\P,2}(n)$, which constitutes Theorem \ref{concgeg}.

\subsection{Urns and marbles}
 Fix $n\geq 2$ an integer. The first thing that one can observe about a hypothetical prime pair $p\leq q$ which sums to $2n$ is that $p\leq n \leq q$. As $p\mid 2n$ implies $p\mid 2n-p$, we have either $(p,2n)=(q,2n)=1$ or $p=q=n$. With this in mind, let
 \begin{gather*}
  \mathcal{A}_n := \{1 < k < n : (k,2n) = 1\}, \\
  \mathcal{B}_n := \{n < k < 2n-1 : (k,2n) = 1\},
 \end{gather*}
 and also
 \[ K(n) := |\mathcal{A}_n|,\qquad P(n) := |\mathcal{A}_n \cap \P|,\qquad Q(n) := |\mathcal{B}_n \cap \P|. \]
 Notice that $k \in \mathcal{A}_n$ if and only if $2n-k \in \mathcal{B}_n$, hence $|\mathcal{A}_n| = |\mathcal{B}_n|$ and we may define the family of bijections:
 \begin{equation*}
  \begin{split}
   \Psi_n: \mathcal{A}_n &\rightarrow \mathcal{B}_n \\
   k~ &\mapsto 2n-k
  \end{split}
 \end{equation*}
 Finally, denote by $g(n):= |\{p\leq n : 2n-p \in \P\}|$ the \emph{Goldbach function}, which satisfies $r_{\P,2}(2n) = 2g(n)-\mathbbm{1}_\P(n)$. This is just a matter of notation, for it will simplify the argument. Rewriting $g(n)$ in terms of $\Psi_n$, we have
 \begin{align*}
  g(n) = \sum_{n\leq p < 2n} \mathbbm{1}_\P(2n-p) &= \Bigg(\sum_{k \in \mathcal{B}_n\cap\P} \mathbbm{1}_\P(2n-k)\Bigg) + \mathbbm{1}_\P(n) \\
  &= \Bigg(\sum_{k \in \mathcal{B}_n\cap\P} \mathbbm{1}_\P(\Psi_n^{-1}(k))\Bigg) + \mathbbm{1}_\P(n) \\
  &= \Bigg(\sum_{k \in \Psi_n^{-1}(\mathcal{B}_n\cap\P)} \mathbbm{1}_\P(k)\Bigg) + \mathbbm{1}_\P(n).
 \end{align*}
 
 The difficulty in understanding the behavior of $g$ lies exactly on the problem of counting primes in the pre-image of $\mathcal{B}_n\cap \P$ by $\Psi_n$. A simple yet instructive attempt to try to grasp this quantity is to consider a random subset of $\mathcal{B}_n$ with $Q(n)$ elements. That is:
 \begin{equation}
  \gtil(n) := \sum_{\substack{k\in \Psi_n^{-1}(X) \\ X\subseteq \mathcal{B}_n\,:\, |X| = Q(n)}} \mathbbm{1}_\P(k). \label{gtil1}
 \end{equation}
 Choosing $X$ uniformly at random, the r.v. $\gtil(n)$ is described by the following parameters: from an urn of $K(n)$ ``marbles'' (set $\mathcal{A}_n$), $P(n)$ are ``special'' (set $\mathcal{A}_n \cap \P$), and $Q(n)$ marbles are drawn uniformly at random and without replacement (set $\Psi_n^{-1}(X)$). The quantity $\gtil(n)$ counts the number of special marbles drawn, and therefore follows a \emph{hypergeometric probability distribution}. This means that
 \begin{equation}
  \Pr(\gtil(n)=k) = \frac{\binom{P(n)}{k}\binom{K(n)-P(n)}{Q(n)-k}}{\binom{K(n)}{Q(n)}}, \label{hprob}
 \end{equation}
 defined for the range
 \[ \max\{0,P(n)+Q(n)-K(n)\} \leq k \leq \min\{P(n),Q(n)\}, \]
 with $\Pr(\gtil(n)=k)=0$ for other values of $k$. The deduction of \eqref{hprob} is rather straightforward, and may be found in Section II.6, p. 43 of Feller \cite{feller1}. It is worth noting that, when $P(n)+Q(n)>K(n)$, we have $\Pr(\gtil(n)=0) = 0$, therefore such $n$ \emph{must}, by the pigeonhole principle, satisfy $g(n)>0$. Writing $K$, $P$, $Q$ in terms of classical arithmetic functions:
 \begin{equation}
  \begin{split}
   K(n) &=\varphi(2n)/2-1, \\
   P(n) &=\pi(n)-\omega(2n), \\
   Q(n) &=\pi(2n-2)-\pi(n),
  \end{split} \label{kpqex}
 \end{equation} 
 one arrives at the following statement.
 
 \begin{prop}\label{gb1}
  If $\pi(2n)-\omega(2n)> \varphi(2n)/2$, then $g(n)>0$.
 \end{prop} 

 Unfortunately, the growth of these arithmetic functions are incompatible with the previous inequality, in the sense that it holds only for finitely many exceptional cases, which cease to occur for $2n$ near $10^5$, as illustrated by the following graph:
 
 \medskip
 \noindent
 \begin{minipage}{\textwidth}
  \centering
   \captionsetup{type=figure}
   \includegraphics[width=.9\textwidth]{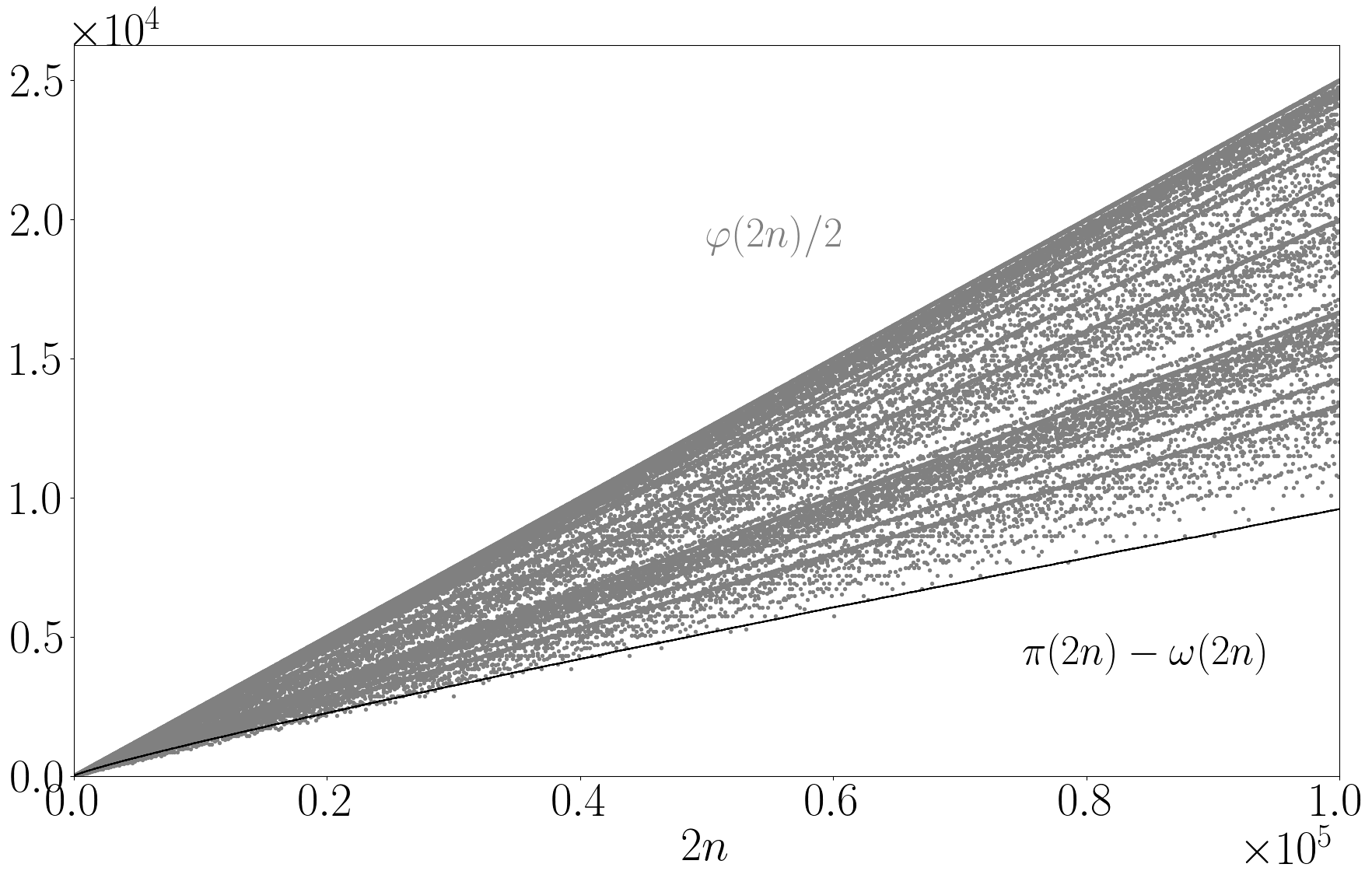}
   \captionof{figure}{Comparison between $\varphi(2n)/2$ and $\pi(2n)-\omega(2n)$ for $2n$ less than $10^5$.}
  \label{fig1}
 \end{minipage}
 \medskip
 
 \noindent
 The last few values of $2n$ for which the hypothesis from Proposition \ref{gb1} holds are $2n=60060$, $78540$ and $90090$. Any improvement to this inequality would potentially require an approach that breaks the symmetric nature of our argument, and for that very reason, this model is incomplete up to this point. The reason is that it neglects primes that do \emph{not} divide $2n$. Amongst the possible ways one could try to correct this,\footnote{In general, the principled way of making such arguments is through sieve theory (cf. Example 5, Section 1.3 of Halberstam and Richert \cite{halberstam11}). For a concrete application in establishing an explicit upper bound to Goldbach's problem, cf. Theorem 3.11, Section 3.7 of \cite{halberstam11}.} we choose the one we deem the simplest, given the informal nature of our argument. For each prime $q\leq \sqrt{2n}$ such that $q\nmid 2n$, if $p\neq q$ is a prime with $2n-p$ also prime, then $p \not\equiv 2n~\mathrm{(mod }~q\mathrm{)}$. Moreover, for any odd prime $p$ we have
 \begin{gather*}
  \lim_{n\to +\infty} \frac{|\{ k\leq n : k \stackrel{p}{\not\equiv} 0 \}|}{n} = 1-\frac{1}{p}, \\
  \lim_{\substack{n\to +\infty \\ p\,\nmid\, n}} \frac{|\{ k\leq n : k \stackrel{p}{\not\equiv} 0 \text{ and } k \stackrel{p}{\not\equiv} 2n \}|}{|\{ k\leq n : k \stackrel{p}{\not\equiv} 0 \}|} = 1-\frac{1}{p-1},
 \end{gather*}
 thus, for every $x\in X$ in \eqref{gtil1}, we may add the correction factor
 \begin{equation}
  ``\Pr\left(\Psi_{n}^{-1}(x) \stackrel{q}{\not\equiv} 0 ~\big|~ x \stackrel{q}{\not\equiv} 0 \right)" \approx \frac{\big(1-\frac{1}{q-1}\big)}{\big(1-\frac{1}{q}\big)} = 1- \frac{1}{(q-1)^2} \label{corrfac}
 \end{equation}
 where $q\leq \sqrt{2n}$ is prime and $q\nmid 2n$. Notice that we wrote $\Pr$ in between quotes, for this is more of a density-type heuristic than a probabilistic one; the more standard version of this argument, explaining the presence of $C_2$ in the twin primes conjecture, can be found in Golomb \cite{golo60}. With this in mind, our \emph{leap of faith} may then be put as follows:
 \begin{equation}
  g(n) \,\stackrel{?}{\approx}\, \gtil(n) \prod_{\substack{p\,\nmid\, 2n \\ p\leq \sqrt{2n}}} \left(1-\frac{1}{(p-1)^2}\right). \label{mainleap}
 \end{equation}
 
 It is important to note that the formula above does not mean anything, for on the left we have an arithmetic function whilst on the right we have a random variable. What we are going to do now is to show that the RHS of this relation is almost surely (i.e. with probability $1$) asymptotic to Hardy \& Littlewood's conjectural estimate for Goldbach's problem. The implied heuristic here, as in most heuristic arguments and conjectures based upon some variation of Cram\'{e}r's model, is that the sequence of primes behaves \emph{pseudorandomly} when analyzed through certain sets of parameters, and hence is expected to fall within the range of certain statistical properties of random sequences.
 
\subsection{Concentration bounds}
 Whenever a r.v. $X$ follows a hypergeometric distribution with parameters $N$ (total marbles in the urn), $M$ (special marbles) and $n$ (drawn marbles) we will say that \emph{$X$ follows $\h(N,M,n)$}, for short.\footnote{This notation is based on M. Skala \cite{skalaHT}.} In our case, from \eqref{gtil1}, we have that $\gtil(n)$ follows $\h(K(n), P(n), Q(n))$. As the expected value for the hypergeometric distribution is just the proportion of special marbles times the size of the sample drawn,\footnote{cf. Example (d) at Section IX.5, pp. 232--233 of Feller \cite{feller1}.} we have
 \begin{equation*}
  \mathbb{E}(\gtil(n)) = \frac{P(n)Q(n)}{K(n)}.
 \end{equation*}
 Going back to \eqref{kpqex}, we have $K(n) = \varphi(2n)/2-1$, and, by the prime number theorem, we have $P(n)$, $Q(n) \sim n\log(n)^{-1}$. Therefore:
 \begin{equation}
  \mathbb{E}(\gtil(n)) \sim \frac{2n}{\varphi(2n)}\frac{n}{\log(n)^2}. \label{expest}
 \end{equation} 
 
 We will now use two tools from probability theory: the already stated Borel--Cantelli lemma, and a concentration inequality for the hypergeometric distribution. The latter is a Hoeffding-type inequality due to V. Chv\'{a}tal \cite{chv79}.
 
 \begin{lem}[Hoeffding--Chv\'{a}tal]\label{hgmt}
  If $X$ is a r.v. that follows $\h(N,M,n)$, then
  \[ \Pr(|X-\mathbb{E}(X)| \geq tn) \leq 2e^{-2t^2n}, \]
  for any real $t\geq 0$.
 \end{lem}
 
 With this, we can state and prove the result at the core of our argument.

 \begin{thm}\label{concgeg}
  As $n\to +\infty$, we have $\gtil(n) \,\uptxt{a.s.}{\sim}\, \mathbb{E}(\gtil(n))$.
 \end{thm}
 \begin{proof}
  Our goal is to apply the Borel--Cantelli lemma, hence, we start by applying Lemma \ref{hgmt} to $\gtil(n)$. First, we have
  \[ \Pr\left(\left|\gtil(n)-\mathbb{E}(\gtil(n))\right| \geq tQ(n) \right) \leq 2e^{-2t^2Q(n)}. \]
  Dividing what is inside by $\mathbb{E}(\gtil(n)) = P(n)Q(n)/K(n)$,
  \[ \Pr\left(\left|\frac{\gtil(n)}{\mathbb{E}(\gtil(n))}-1\right| \geq \frac{tK(n)}{P(n)} \right) \leq 2e^{-2t^2Q(n)}. \]
  From \eqref{kpqex}, we see that $K(n)< n$ for all $n\geq 2$, and from the PNT we have $P(n)$, $Q(n) > n/2\log(n)$ for all sufficiently large $n$. Therefore,
  \[ \Pr\left(\left|\frac{\gtil(n)}{\mathbb{E}(\gtil(n))}-1\right| \geq \frac{nt}{n/2\log(n)} \right) \ll e^{-nt^2/2\log(n)}. \]
  For each $n\geq 2$, change $t$ to $1/2\log(n)^2$. Thus, we obtain
  \[ \Pr\left(\left|\frac{\gtil(n)}{\mathbb{E}(\gtil(n))}-1\right| \geq \frac{1}{\log(n)} \right) \ll e^{-n/8\log(n)^5}. \]
  Since $\log(n)^{-1} = o(1)$ and $\sum_{n=2}^{\infty} e^{-n/8\log(n)^5}$ converges, it follows, for all $\varepsilon > 0$,
  \[ \sum_{n\geq 3} \Pr\left(\left|\frac{\gtil(n)}{\mathbb{E}(\gtil(n))}-1\right| \geq \varepsilon \right) < +\infty, \]
  which, by the Borel--Cantelli lemma, implies our theorem.
 \end{proof}

 Hence, in view of \eqref{mainleap}, we are naturally led to expect the following:
 
 \begin{conj}[Heuristic conclusion]\label{heurs}
  As $n\to +\infty$,
  \begin{equation}
   g(n) \sim \frac{2n}{\varphi(2n)}\frac{n}{\log(n)^2} \prod_{\substack{p\,\nmid\, 2n \\ p \leq \sqrt{2n}}} \left(1-\frac{1}{(p-1)^2}\right). \label{heursEQ}
  \end{equation}
 \end{conj}
 
 The term on the RHS is just the product of \eqref{corrfac} and \eqref{expest}. From the fact that $r_{\P,2}(2n) = 2g(n)-\mathbbm{1}_{\P}(n)$, the estimate at \eqref{heursEQ} may be rewritten in the notation used in our introductory remarks as follows: \emph{when $n\to+\infty$ through the even numbers,}
 \begin{equation}
  r_{\P,2}(n) \sim \frac{n}{\varphi(n)}\frac{n}{\log(n)^2} \prod_{\substack{p\,\nmid\, n \\ 3\leq p \leq \sqrt{n}}} \left(1-\frac{1}{(p-1)^2}\right). \label{gconj}
 \end{equation}
 Notice, however, that when $n$ is even,
 \begin{equation*}
  \frac{n}{\varphi(n)} = \frac{n}{n\prod\limits_{p\mid n}\big(1-\frac{1}{p}\big)}
  = 2\prod_{\substack{p\mid n \\ p\geq 3}} \left(\frac{p}{p-1}\right) = 2 \prod_{\substack{p\mid n \\ p\geq 3}} \left(1-\frac{1}{(p-1)^2}\right) \prod_{\substack{p\mid n \\ p\geq 3}} \left(\frac{p-1}{p-2}\right),
 \end{equation*}
 thus, \eqref{gconj} may be rewritten as
 \begin{equation}
  r_{\P,2}(n) \sim 2\prod_{3\leq p \leq \sqrt{n}} \left(1-\frac{1}{(p-1)^2}\right)\frac{n}{\log(n)^2} \prod_{\substack{p\mid n \\ p\geq 3}} \left(\frac{p-1}{p-2}\right), \label{ggconj}
 \end{equation}
 which, by the convergence of $C_2$ (as in \eqref{c2}), is equivalent to Hardy \& Littlewood's conjectured estimate \eqref{hlgbcj}. To finish our discussion, we present two graphs showing how accurate this estimate is for small values of $n$. In our second graph, $\int_{2}^{n-2} (\log(t)\log(n-t))^{-1}\mathrm{d}t$ is used instead of $n\log(n)^{-2}$, for it gives a faster approximation. (This does not change the conclusion, for $\int_{2}^{n-2} \frac{\mathrm{d}t}{\log(t)\log(n-t)} \sim n\log(n)^{-2}$.)
 
 { \newpage\null\vfill
 \smallskip
 \noindent
 \begin{minipage}{\textwidth}
  \centering
   \captionsetup{type=figure}
   \includegraphics[width=\textwidth]{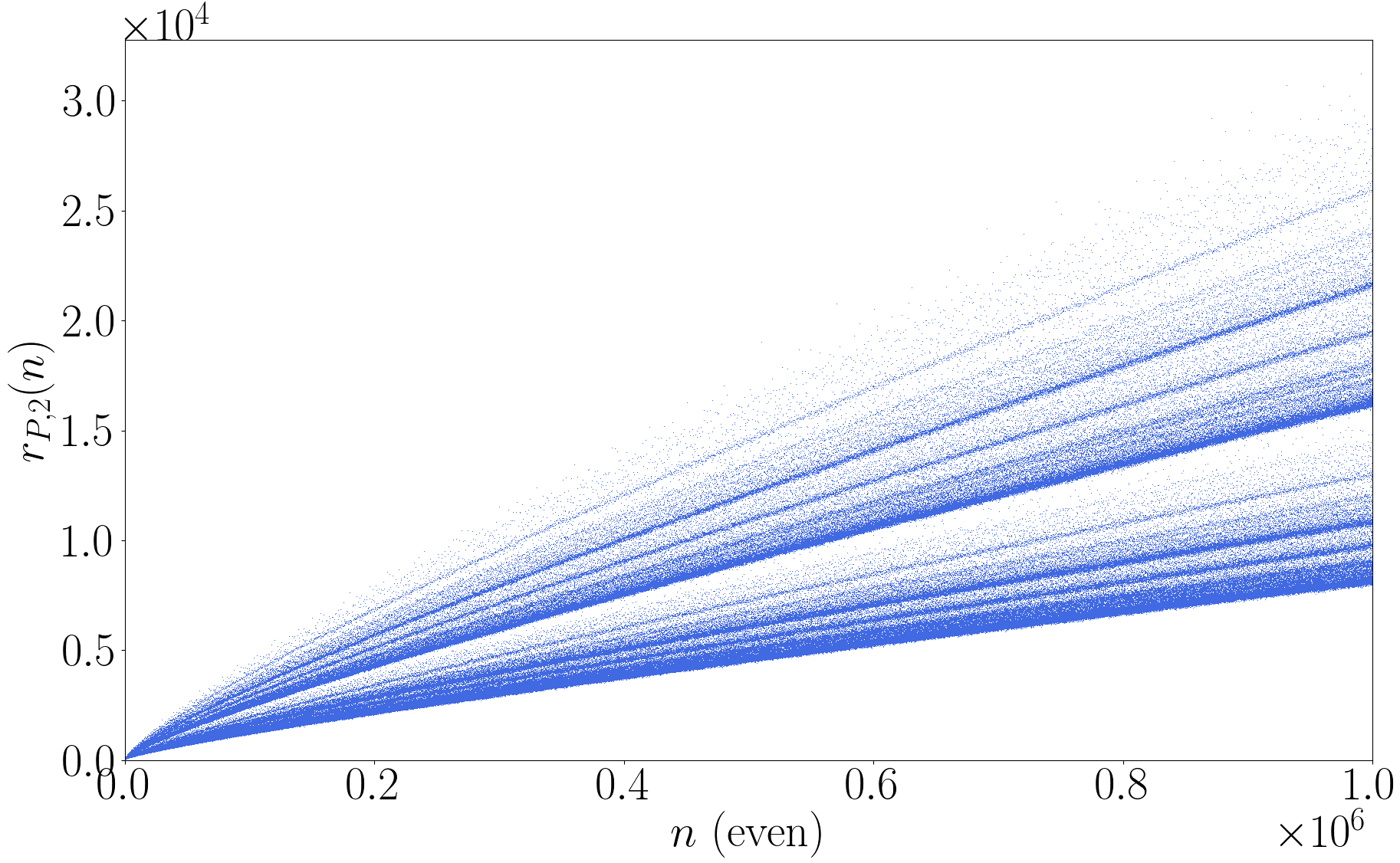}
   \captionof{figure}{Values of $r_{\P,2}(n)$ for $n=2$ to $10^6$ (even).}
   \label{fig3}
 \end{minipage}
 
 \vfill
 
 \noindent
 \begin{minipage}{\textwidth}
  \centering
   \captionsetup{type=figure}
   \includegraphics[width=\textwidth]{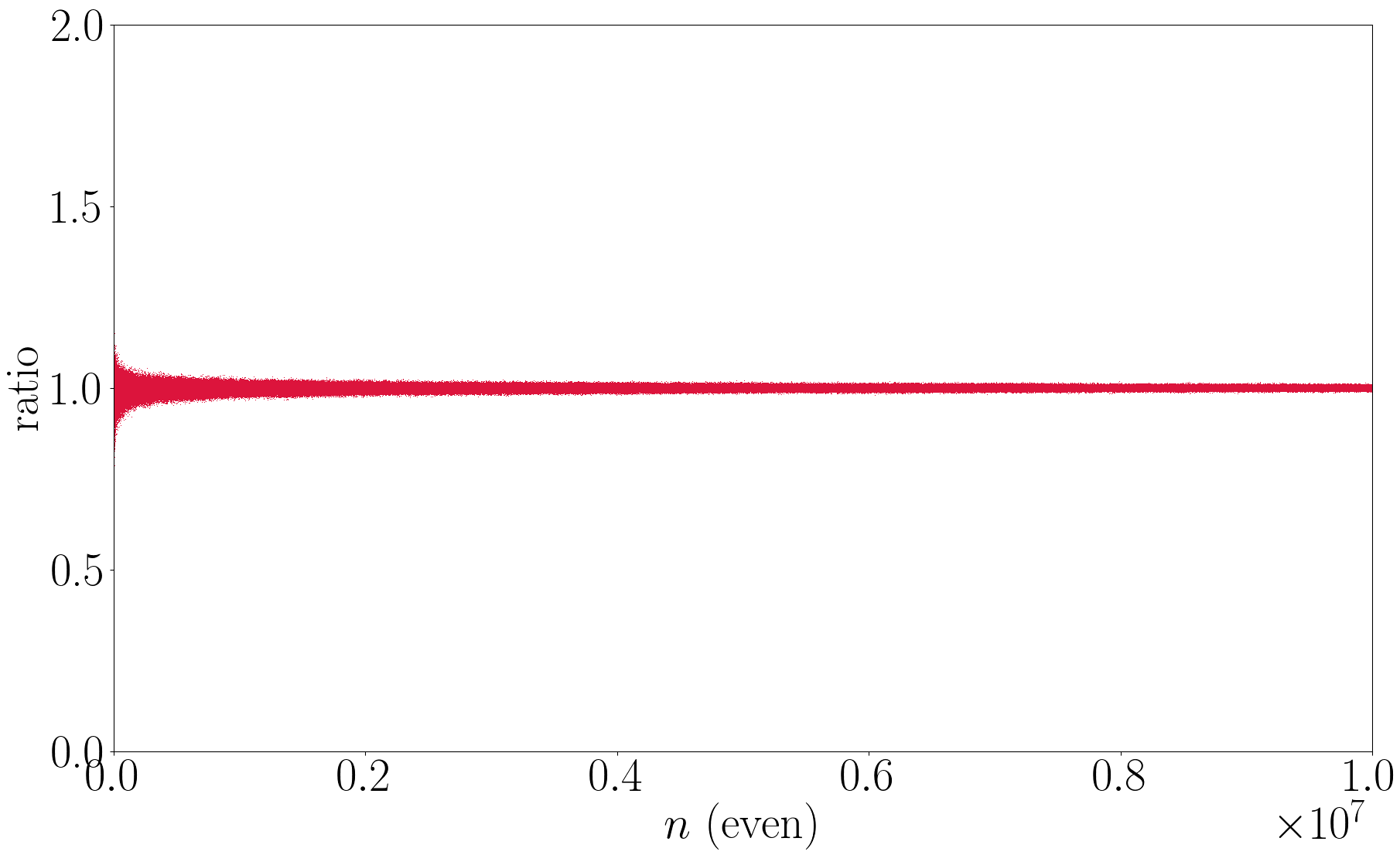}
   \captionof{figure}{From $n=2$ to $10^7$ (even), ratio between $r_{\P,2}(n)$ and $2 C_2 \big(\int_{2}^{n-2} \frac{\mathrm{d}t}{\log(t)\log(n-t)}\big) \prod_{\underset{p\geq 3}{p\mid n}} \big(\frac{p-1}{p-2}\big)$.}
   \label{fig2}
 \end{minipage}
 \vfill\newpage }

\bibliographystyle{amsplain}
\bibliography{$HOME/Acad/Writings/_latex/bibliotheca}%

\providecommand{\bysame}{\leavevmode\hbox to3em{\hrulefill}\thinspace}
\providecommand{\MR}{\relax\ifhmode\unskip\space\fi MR }
\providecommand{\MRhref}[2]{%
  \href{http://www.ams.org/mathscinet-getitem?mr=#1}{#2}
}
\providecommand{\href}[2]{#2}
\begin{thebibliography}{10}

\bibitem{apostol76}
T.~M. Apostol, \emph{An introduction to analytic number theory}, Springer,
  1976.

\bibitem{bathor62}
P.~T. Bateman and R.~A. Horn, \emph{A heuristic asymptotic formula concerning
  the distribution of prime numbers}, Math. Comput. \textbf{16} (1962),
  363--367.

\bibitem{chv79}
V.~Chv\'{a}tal, \emph{The tail of the hypergeometric distribution}, Discrete
  Math. \textbf{25} (1979), 285--287.

\bibitem{cra36}
H.~Cram\'{e}r, \emph{On the order of magnitude of the difference between
  consecutive prime numbers}, Acta Arith. \textbf{2} (1936), 23--46.

\bibitem{deshenlan98}
J.-M. Deshouillers, F.~Hennecart, and B.~Landreau, \emph{Sums of powers: an
  arithmetic refinement to the probabilistic model of {Erd\H{o}s} and
  {R\'{e}nyi}}, Acta Arith. \textbf{85} (1998), 13--33.

\bibitem{erdo56}
P.~Erd\H{o}s, \emph{Problems and results in additive number theory}, Colloque
  sur la Theorie des Nombres (CBRM) (Bruxelles), 1956, pp.~127--137.

\bibitem{erdren60}
P.~Erd\H{o}s and A.~R\'{e}nyi, \emph{Additive properties of random sequences of
  positive integers}, Acta Arith. \textbf{6} (1960), 83--110.

\bibitem{erdtet90}
P.~Erd\H{o}s and P.~Tetali, \emph{Representations of integers as the sum of $k$
  terms}, Random Struct. Algor. \textbf{1} (1990), 245--261.

\bibitem{feller1}
W.~Feller, \emph{An introduction to probability theory and its applications},
  3rd ed., vol.~1, John Wiley \& Sons Inc., New York, 1967.

\bibitem{golo60}
S.~W. Golomb, \emph{The twin prime constant}, Amer. Math. Monthly \textbf{67}
  (1960), no.~8, 767--769.

\bibitem{gra95}
A.~Granville, \emph{{Harald Cram\'{e}r} and the distribution of prime numbers},
  Scand. Actuar. J. \textbf{1} (1995), 12--28.

\bibitem{guy94}
R.~K. Guy, \emph{Unsolved problems in number theory}, 2nd ed., Springer, 1994.

\bibitem{halberstam11}
H.~Halberstam and H.-E. Richert, \emph{Sieve methods}, dover ed., Dover
  Publications, New York, 2011.

\bibitem{halberstam83}
H.~Halberstam and K.~F. Roth, \emph{Sequences}, revised ed., Springer, 1983.

\bibitem{halmos74}
P.~R. Halmos, \emph{Measure theory}, Graduate Texts in Mathematics, vol.~18,
  Springer, 1974.

\bibitem{harlit23}
G.~H. Hardy and J.~E. Littlewood, \emph{Some problems of 'partitio numerorum'
  iii: on the expression of a number as a sum of primes}, Acta Mathematica
  \textbf{44} (1923), 1--70.

\bibitem{harlit24}
\bysame, \emph{Some problems of 'partitio numerorum' (v): A further
  contribution to the study of {G}oldbach's problem}, P. Lond. Math. Soc.
  \textbf{22} (1924), no.~2, 46--56.

\bibitem{hardy08}
G.~H. Hardy and E.~M. Wright, \emph{An introduction to the theory of numbers},
  6th ed., Oxford University Press, 2008.

\bibitem{hua65}
L.~K. Hua, \emph{Additive theory of prime numbers}, American Mathematical
  Society, 1965.

\bibitem{lan95}
B.~Landreau, \emph{\'{E}tude probabiliste des sommes des puissances
  $s$-i\`{e}mes}, Compos. Math. \textbf{99} (1995), 1--31.

\bibitem{monvau75}
H.~L. Montgomery and R.~C. Vaughan, \emph{The exceptional set in {G}oldbach's
  problem}, Acta Arith. \textbf{27} (1975), 353--370.

\bibitem{nathanson96}
M.~B. Nathanson, \emph{Additive number theory: The classical bases}, 2nd ed.,
  Graduate Texts in Mathematics, vol. 164, Springer, 1996.

\bibitem{pin07}
J.~Pintz, \emph{{C}ram\'{e}r vs. {C}ram\'{e}r. {O}n {C}ram\'{e}r's
  probabilistic model for primes}, Funct. Approx. Comment. Math. \textbf{37}
  (2007), 361--376.

\bibitem{skalaHT}
M.~Skala, \emph{Hypergeometric tail inequalities: ending the insanity},
  \href{https://arxiv.org/abs/1311.5939v1}{\texttt{arXiv:1311.5939v1}}
  [math.PR], 2013.

\bibitem{taf19}
C.~T\'{a}fula, \emph{An extension of the {Erd\H{o}s--Tetali} theorem}, Random
  Struct. Algor. \textbf{55} (2019), no.~1, 173--214.

\bibitem{tao06}
T.~Tao and V.~H. Vu, \emph{Additive combinatorics}, Cambridge Stud. Adv. Math.,
  vol. 105, Cambridge Univ. Press, 2006.

\bibitem{vinogradov54}
I.~M. Vinogradov, \emph{The method of trigonometrical sums in the theory of
  numbers}, Interscience, London, 1954, translated by {K. F. Roth} and {A.
  Davenport}.

\bibitem{vvu00wp}
V.~H. Vu, \emph{On a refinement of {W}aring's problem}, Duke Math. J.
  \textbf{105} (2000), 107--134.

\end{thebibliography}
\end{document}